\input amstex\documentstyle{amsppt}  
\pagewidth{12.5cm}\pageheight{19cm}\magnification\magstep1
\topmatter
\title From families in Weyl groups to Springer representations\endtitle
\author G. Lusztig\endauthor
\address{Department of Mathematics, M.I.T., Cambridge, MA 02139}\endaddress
\thanks{Supported by NSF grant DMS-1855773.}\endthanks
\endtopmatter   
\document

\define\Irr{\text{\rm Irr}}

\define\mpb{\medpagebreak}

\define\si{\sim}

\define\sqc{\sqcup}

\define\part{\partial}
\define\emp{\emptyset}

\define\m{\mapsto}
\define\do{\dots}

\define\sm{\smallmatrix}
\define\esm{\endsmallmatrix}
\define\sub{\subset}    
\define\bxt{\boxtimes}
\define\T{\times}
\define\ti{\tilde}
\define\nl{\newline}

\define\ot{\otimes}
\define\bbq{\bar{\QQ}_l}

\define\Hom{\text{\rm Hom}}

\define\Ind{\text{\rm Ind}}

\define\g{\gamma}
\redefine\d{\delta}
\define\e{\epsilon}

\redefine\G{\Gamma}

\define\ee{\bold e}

\define\kk{\bold k}

\define\CC{\bold C}

\define\NN{\bold N}

\define\QQ{\bold Q}

\define\cb{\Cal B}

\define\cj{\Cal J}

\define\cx{\Cal X}

\define\fc{\frak c}

\define\sha{\sharp}

\define\bul{\bullet}

\define\che{\check}

\head 1. The function $c:\Irr W@>>>\NN$\endhead
\subhead 1.1\endsubhead
Let $G$ be an almost simple, simply connected algebraic group over an algebraically closed field $\kk$
of characteristic $p\ge0$. Let $W$ be the Weyl group of $G$.
In this section and until the end of 3.3 we assume that $p$ is $0$ or a good prime for $G$.

For any unipotent class $C$ of $G$ let $E_C$
be the representation of $W$ defined by Springer \cite{S76} on the top $l$-adic cohomology
of the Springer fibre $\cb_u$ at an element $u\in C$. (Here $l$ is fixed prime number $\ne p$.)
For any finite group $\G$ let $\Irr \G$ be a set of representatives for the irreducible representations
of $\G$
over the $l$-adic numbers and let $R_{\G}$ be the free abelian group with basis $\Irr\G$. We can view $E_C$
as an element of $R_W$. According to Springer, there is a unique partition $\Irr W=\sqc_C\Irr_CW$ (the
``Springer partition'') where $C$ runs over the unipotent classes of $G$ such that $\Irr_CW$ consists of all
$E\in\Irr W$ which appear in $E_C$ with coefficient $>0$. We define $\g:\Irr W@>>>\NN$ by $\g(E)=\dim\cb_u$
where $E\in\Irr_CW$ and $u\in C$.
The purpose of this paper is to propose a definition of the Springer partition, of the
collection of representations $E_C$ and of the function $\g:\Irr W@>>>\NN$ which is purely algebraic
(without use of geometry) and which is suitable for computer calculations. The same objects
can be obtained (without use of geometry) from the approach of \cite{LY}; but that approach
is more complicated than the present one and seems to be unsuitable for computer calculations. 

\subhead 1.2\endsubhead
Let $W'$ be a Weyl group. For $E'\in\Irr W'$
the generic degree of the irreducible Hecke algebra representation $E'_q$ corresponding to
$E'$ is of the form $(1/m_{E'})q^{a_{E'}}+$ terms of strictly higher degree in $q$. (Here $m_{E'}$ is a
constant integer $\ge1$.)
In \cite{L82} we have defined a partition of $\Irr W'$ into
subsets called {\it families}. The definition is purely algebraic; it involves induction from parabolic
subgroups, tensoring by the sign representation and the knowledge of the function $\Irr W'@>>>\NN$,
$E'\m a_{E'}$. It is known that $E'\m a_{E'}$ is constant on each family. 
Let $\Irr_{sp}W'$ be the subset of $\Irr W'$ consisting of the special
representations. If $E'\in\Irr W'$ then there is a unique $E'_0\in\Irr_{sp}W'$ in the
same family as $E'$ and $n_{E'}=m_{E'_0}/m_{E'}$ is an integer $\ge1$.

\subhead 1.3\endsubhead
Let $W_{af}$ be the affine Weyl group associated to the reductive group over $\CC$ of type dual to that of $G$.
Let $S_{af}$ be the set of simple reflections of $W_{af}$. We can identify $W$ with the quotient of $W_{af}$
by the group of translations. For any $K\subsetneqq S_{af}$, the subgroup of $W_{af}$ generated by $K$ is a
finite Weyl group; we will identify it with its image $W_K$ under the canonical surjection $W_{af}@>>>W$.

For $E\in\Irr W$ let $\G(E)$ be the set of all pairs $(K,E')$ where $K\subsetneqq S_{af}$ and $E'\in\Irr W_K$
is such that the multiplicity $(E:\Ind_{W_K}^W(E'))$ is nonzero. We have $\G(E)\ne\emp$. Let
$c_E=\max\{a_{E'};(K,E')\in\G(E)\}$; here $a_{E'}$ is computed in terms of $W_K$. Let
$\G_M(E)=\{(K,E')\in\G(E);a_{E'}=c_E\}$. We have $\G_M(E)\ne\emp$. Let $N_E=\max\{n_{E'};(K,E')\in\G_M(E)\}$.
Let $\si$ be the equivalence relation on $\Irr W$ generated by the relation $E_1\si E_2$ if there exist
$(K_1,E'_1)\in\G_M(E_1),(K_2,E'_2)\in\G_M(E_2)$ such that $K_1=K_2$ and $E'_1,E'_2$ are in the same family of
$W_{K_1}=W_{K_2}$. The equivalence classes on $\Irr W$ for $\si$ are called $c$-families. Note that the
function $E\m c_E$ is constant on $c$-families. For any $c$-family $\frak f\sub\Irr W$ we set
$\ee_{\frak f}=\sum_{E\in\frak f}N_EE\in R_W$.

In the definition of $\G(E)$ one could add the condition that $\sha(K)=\sha(S_{af})-1$. This would not
affect the validity of the following result.

\proclaim{Proposition 1.4}(a) For any $E\in\Irr W$ we have $c_E=\g(E)$.

(b) The Springer partition of $\Irr W$ coincides with the partition of $\Irr W$ into $c$-families.

(c) The collection of Springer representations $E_C$ of $W$ for various unipotent classes $C$ coincides with
the collection of elements $\ee_{\frak f}\in R_W$ for various $c$-families $\frak f$  in $\Irr W$.
\endproclaim

\subhead 1.5\endsubhead
For any $K\subsetneqq S_{af}$ we define a homomorphism $\cj^W_{W_K}:R_{W_K}@>>>R_W$ by setting for 
any $E'\in\Irr W_K$
$$\cj^W_{W_K}(E')=\sum_{E\in\Irr W;c_E=a_{E'}}(E:\Ind_{W_K}^W(E'))E\in R_W.$$
Let $\phi$ be a family of $W_K$ and let $\frak f$ be a $c$-family. We say that $(K,\phi)$ is {\it adapted} to
$\frak f$ if there exists a subset $\phi_*$ of $\phi$ such that

(a) $\cj^W_{W_K}(E')\in\frak f$ (and in particular is in $\Irr W$) for any $E'\in\phi_*$;

(b) $E'\m E:=\cj^W_{W_K}(E')$ is a bijection $\phi_*@>\si>>\frak f$;

(c) if $E',E$ are as in (b) then $n_{E'}=N_E$;

(d) if $E'\in\phi-\phi_*$ and $E\in\Irr W$ appears in $\Ind_{W_K}^W(E')$, then $c_E>a_{E'}$.
\nl
Note that $\phi_*$ above is unique (if it exists). Moreover the value of the $a$-function on $\phi$
must be equal to the value of the $c$-function on $\frak f$.
For $(K,\phi)$ adapted to $\frak f$ we set $E'_{\phi_*}=\sum_{E'\in\phi_*}n_{E'}E'\in R_{W_K}$.
We then have $\cj^W_{W_K}(E'_{\phi_*})=\ee_{\frak f}$. Note that $\phi_*$ is not determined by the
Coxeter group $W_K$ (it depends on $\frak f$).

\proclaim{Proposition 1.6} Let $\frak f\sub\Irr W$ be a $c$-family. There exist $K\subsetneqq S_{af}$ and a
family $\phi$ of $W_K$ such that $(K,\phi)$ is adapted to $\frak f$.
\endproclaim

\subhead 1.7 \endsubhead
Let $W'$ be a Weyl group; let $S'$ be the set of simple reflections in $W'$. For any $H\sub S'$ let
$W'_H$ be the subgroup of $W'$ generated by $H$. Let $E'\in\Irr_{sp}W'$. We say that $E'$ is 
non-rigid if there exists $H\subsetneqq S'$ and $E''\in\Irr_{sp}W'_H$ such that
$E'$ appears in $\Ind_{W_H}^W(E'')$ and $a_{E''}=a_{E'}$, $m_{E''}=m_{E'}$. (Here $a_{E''},m_{E''}$ are
defined relative to $W'_H$.) We say that $E'\in\Irr_{sp}W'$ is rigid if it is not non-rigid.
For any $E'\in\Irr_{sp}W'$ one can find $H\sub S'$ and $E''\in\Irr_{sp}W'_H$ such that
$E''$ is rigid, $E'$ appears in $\Ind_{W_H}^W(E'')$ and $a_{E''}=a_{E'}$, $m_{E''}=m_{E'}$.
Moreover, $(H,E'')$ is uniquely determined by $E'$ up to conjugation by an element of $W'$.
(A statement close to this appears in \cite{L84, 13.1}.)
A family of $W'$ is said to be rigid if the unique special representation in it is rigid.
If $W'$ is a product of two Weyl groups $W'_1,W'_2$, the rigid special representations of $W'$ are
precisely the external tensor products of a rigid special representation of $W'_1$ with one of $W'_2$.

We denote by $\e$ the sign representation of $W'$. It is rigid special. If $W'$ is of type $A_n$, $\e$
is the only rigid special representation of $W'$.

In the next subsection we list the rigid special representation for $W'$ irreducible of low rank.
We use the following notation.
If $W'$ is of type $E_6,E_7$ or $E_8$, an element $E'\in\Irr W'$ can be represented uniquely in
the form $\d_x$ where $d=\dim(E')$ and $x$ is the smallest integer $\ge0$ such that $E'$ appears in
the $x$-th symmetric power of the reflection representation of $W'$. The same notation can be used
in type $F_4,G_2$ but in these cases there may be two elements $E'\in\Irr W'$ with the same
$\d_x$. (This ambiguity does not appear for rigid special representations.)
If $W'$ is of type $D_n,n\ge4$ we represent an $E'\in\Irr W'$ as a sequence 
$a_1a_2...a_{2k}$ where $[a_{2k}<a_{2k-2}<\do<a_2,a_{2k-1}<a_{2k-3}<\do<a_1]$ is the
symbol representing $E'$ in \cite{L84}.
If $W'$ is of type $B_n,n\ge2$ we represent an $E'\in\Irr W'$ as a sequence 
$a_0a_1a_2...a_{2k}$ where $[a_{2k}<a_{2k-2}<\do<a_0,a_{2k-1}<a_{2k-3}<\do<a_1]$ is the
symbol representing $E'$ in \cite{L84}. When $n=2$ this notation depends on the order of $S'$;
the representations $201,120$ are interchanged when the order of $S'$ is reversed; again this
ambiguity does not appear for rigid special representations.

\subhead 1.8 \endsubhead
Here are the rigid special representations for $W'$ assumed to be irreducible of low rank.

Type $B_2$: $210$, $22110=\e$ (with $a=1,4$ respectively);

Type $B_3$: $32110$, $3322110=\e$  (with $a=4,9$ respectively);

Type $B_4$: $32210$, $4322110$, $443322110=\e$ (with $a=6,9,16$ respectively);

Type $D_4$:  $3210$, $332110$, $43322110=\e$ (with $a=3,7,12$ respectively);

Type $D_5$:  $432110$, $44322110$, $5443322110=\e$ (with $a=7,13,20$ respectively);

Type $D_6$:  $432210$, $54322110$, $44332110$, $5543322110$, $655443322110=\e$

(with $a=10,13,16,21,30$ respectively);

Type $D_7$:  $433210$, $54332110$, $6543322110$, $5544322110$, $665443322110$,

$76655443322110=\e$ (with $a=12,16,21,24,29,42$ respectively);

Type $D_8$:  $443210$, $54432110$, $54332210$, $6544322110$, $5544332110$, $765443322110$,

$665543322110$, $6655443322110$, $77655443322110$, $8776655443322110=\e$

(with $a=13,18,21,24,29,31,34,42,43,56$ respectively);

Type $E_6$: $80_7,30_{15},6_{25},1_{36}=\e$;

Type $E_7$: $512_{11},315_{16},120_{25},56_{30},27_{37},7_{46},1_{63}=\e$;

Type $E_8$: $4480_{16},4200_{24},4096_{26},2240_{28},1400_{32},1400_{37},700_{42},560_{47},210_{52}$,

$112_{63},35_{74},8_{91},1_{120}=\e$;

Type $F_4$: $12_4,9_{10},4_{13},1_{24}=\e$;

Type $G_2$:  $2_1,1_6=\e$.

\mpb

We now state a refinement of Proposition 1.6.

\proclaim{Proposition 1.9} Let $\frak f\sub\Irr W$ be a $c$-family. There exist
$K\subsetneqq S_{af}$ and a rigid family
$\phi$ of $W_K$ such that $(K,\phi)$ is adapted to $\frak f$.
\endproclaim

\head 2. Tables\endhead
\subhead 2.1\endsubhead
In the tables below, for $W$ of type $E_8,E_7,E_6,F_4,G_2$ we make a list of 
of the elements $\ee_{\frak f}$ (see 1.3) for each $c$-family $\frak f$ of $W$. (Note that
$\frak f$ is determined by $\ee_{\frak f}$.) In each case we five the value of the $c$-function on
$\frak f$ and we specify the pairs $(K,\phi)$ with $\phi$ a rigid family of $W_K$ such
that $(K,\phi)$ is adapted to $\frak f$.
When there are several such pairs $(K,\phi)$ which are obtained one from another by conjugating by
an element of $W$ or by applying an automorphism of the affine diagram of $W_{af}$, we list only
one pair out of these. In each case there remain one or two pairs to be listed.
In each case we specify $K$ by the type of the corresponding subdiagram, as a name or as a picture
(with the elements of $K$ marked by $\bul$). The irreducible representations which enter in $\phi$
are denoted using the conventions in 1.7.
In type $F_4$ the notation $\d_x$ in 1.7 is ambiguous; when it represents two objects, we distinguish 
them by denoting them by $\d'_x,\d''_x$ (in a way compatible with the conventions in
\cite{Ch}); the table in 2.5 could be taken as definition of $d'_x,d''_x$ as well as a definition
of $201,120$ for type $B_2$.
In type $G_2$ there are two objects represented by $1_3$; we write them as
$1_3,\ti1_3$ (they are defined by the table in 2.6).

\subhead 2.2. Table for type $E_8$\endsubhead

$c=0$;  $1_0=\cj_{W_K}^W(\e),K=\emp$;

$c=1$;   $8_1=\cj_{W_K}^W(\e),K=\left(\sm\bul&.&.&.&.&.&.&.\\{}&{}&.&{}&{}&{}&{}&{}\esm\right)$;

$c=2$;   $35_2=\cj_{W_K}^W(\e),
K=\left(\sm.&\bul&.&\bul&.&.&.&.\\{}&{}&.&{}&{}&{}&{}&{}\esm\right)$;

$c=3$;   $112_3+28_8=\cj_{W_K}^W(3210+3201)$,
$K=\left(\sm.&\bul&\bul&\bul&.&.&.&.\\{}&{}&\bul&{}&{}&{}&{}&{}\esm\right)$;

$c=4$;   $84_4=\cj^W_{W_K}(\e),K=\left(\sm.&.&.&\bul&.&\bul&.&\bul\\{}&{}&\bul&{}&{}&{}&{}&{}\esm\right)$;
   
$c=4$;  $210_4+160_7=\cj^W_{W_K}(3210\bxt\e+3201\bxt\e)$,
$K=\left(\sm\bul&\bul&\bul&\bul&.&\bul&.&.\\{}&{}&\bul&{}&{}&{}&{}&{}\esm\right)$;

 $c=5$;  $560_5+50_8=\cj^W_{W_K}(3210\bxt\e+3201\bxt\e)$, $K=D_4\T(A_1\T A_1)$;

$c=6$;   $567_6=\cj^W_{W_K}(\e),K=\left(\sm\bul&\bul&\bul&.&.&.&.&.\\{}&{}&.&{}&{}&{}&{}&{}\esm\right)$;

$c=6$;  $700_6+300_8=\cj^W_{W_K}(3210\bxt\e+3201\bxt\e), K=D_4\T A_2$;

$c=7$;  $400_7=\cj^W_{W_K}(\e),
K=\left(\sm\bul&\bul&.&\bul&.&\bul&.&\bul\\{}&{}&\bul&{}&{}&{}&{}&{}\esm\right)$;
               
$c=7$;  $1400_7+2\T1008_9+56_{19}=\cj^W_{W_K}(80_7+2\T90_8+20_{10}), K=E_6$;

$c=8$; $1344_8=\cj^W_{W_K}(\e),K=\left(\sm.&.&.&\bul&.&\bul&\bul&\bul\\{}&{}&\bul&{}&{}&{}&{}&{}\esm\right)$;
                         
$c= 8$;  $1400_8+2\T1575_{10}+350_{14}=\cj^W_{W_K}(80_7\bxt\e+2\T90_8\bxt\e+20_{10}\bxt\e)$,

$K=E_6\T A_1$;

$c=9$;   $3240_9+1050_{10}=\cj^W_{W_K}(3210\bxt\e+3120\bxt\e)$

$=\cj^W_{W_{K'}}(432110\bxt\e+423110\bxt\e),K=D_4\T A_3, K'=D_5\T(A_1\T A_1)$;

$c=9$;    $448_9=\cj^W_{W_K}(\e),K=\left(\sm\bul&\bul&.&\bul&\bul&.&\bul&\bul
\\{}&{}&.&{}&{}&{}&{}&{}\esm\right)$;

$c=10$;   $2240_{10}+2\T175_{12}+840_{13}=\cj^W_{W_K}(80_7\bxt\e+2\T10_9\bxt\e+20_{10}\bxt\e),K=E_6A_2$;

$c=10$;  $2268_{10}+1296_{13}=\cj^W_{W_K}(432210+432201),K=D_6$;

$c=11$;   $4096_{11}+4096_{12}=\cj^W_{W_K}(512_{11}+512_{12})
=\cj^W_{W_{K'}}(432210\bxt\e+432201\bxt\e)$,

$K=E_7,K'=D_6A_1$;

$c=11$;   $1400_{11}=\cj^W_{W_K}(\e),
K=\left(\sm\bul&\bul&.&\bul&.&\bul&\bul&\bul\\{}&{}&\bul&{}&{}&{}&{}&{}\esm\right)$;

$c=12$;   $525_{12}=\cj^K_{W_K}(\e), K=D_4$;

$c=12$;   $972_{12}=\cj^K_{W_K}(\e), K=\left(\sm\bul&\bul&\bul&.&\bul&\bul&\bul&.
\\{}&{}&.&{}&{}&{}&{}&{}\esm\right)$;

$c=12$;  $4200_{12}+3360_{13}=\cj^W_{W_K}(433210+433201)=\cj^W_{W_{K'}}(512_{11}\bxt\e+512_{12}\bxt\e)$,

$K=D_7,K'=E_7\T A_1$;

$c=13$;  $4536_{13}+840_{14}=\cj^W_{W_K}(432110\bxt\e+431201\bxt\e)
=\cj^W_{W_{K'}}(443210+443201)$,

$K=D_5\T A_3,K'=D_8$;

$c=13$;  $2800_{13}+2100_{16}=\cj^W_{W_K}(54322110+54231201), K=D_6$;

$c=14$;  $6075_{14}+700_{16}=\cj^W_{W_K}(54322110\bxt\e+53422110\bxt\e),K=D_6\T A_1$;

$c=14$;  $2835_{14}=\cj^W_{W_K}(\e),
K=\left(\sm\bul&.&\bul&.&\bul&\bul&\bul&\bul\\{}&{}&\bul&{}&{}&{}&{}&{}\esm\right)$;

$c=15$;  $4200_{15}=\cj^W_{W_K}(\e),
K=\left(\sm.&\bul&\bul&\bul&.&\bul&\bul&.\\{}&{}&\bul&{}&{}&{}&{}&{}\esm\right)$;

$c=15$;  $5600_{15}+2400_{17}=\cj^W_{W_K}(30_{15}+15_{17}), K=E_6$;

$c=16$;  $4480_{16}+5\T4536_{18}+4\T5670_{18}+5\T1400_{20}+6\T1680_{22}+4\T70_{32}
=\cj^W_{W_K}(4480_{16}+5\T4536_{18}+4\T5670_{18}+5\T1400_{20}+6\T1680_{22}+4\T70_{32}),K=E_8$;

$c=16$;  $3200_{16}=c^W_{W_K}(\e), K=A_5\T A_1$;

$c=17$;   $7168_{17}+5600_{19}+448_{25}=\cj^W_{W_K}(315_{16}\bxt\e+2\T280_{18}\bxt\e+35_{22}\bxt\e)$,

$K=E_7\T A_1$;

$c=18$;  $3150_{18}+1134_{20}=\cj^W_{W_K}(30_{15}\bxt\e+15_{17}\bxt\e),K=E_6\T A_2$;

$c=18$;  $4200_{18}+2688_{20}=\cj^W_{W_K}(54432110+54431201), K=D_8$;

$c=19$;  $1344_{19}=\cj^W_{W_K}(44322110\bxt\e), K=D_5\T A_3$;

$c=19$;  $2016_{19}=\cj^W_{W_K}(\e), K=A_5\T A_2\T A_1$;

$c=20$;   $420_{20}=\cj^W_{W_K}(\e),K=A_4\T A_4$;

$c=20$;   $2100_{20}=\cj^W_{W_K}(\e), K=\left(\sm\bul&\bul&\bul&\bul&.&.&.&.
\\{}&{}&\bul&{}&{}&{}&{}&{}\esm\right)$;

$c=21$;   $4200_{21}+168_{24}=\cj^W_{W_K}(54332210+53423120), K=D_8$;

$c=21$;  $5600_{21}+2400_{23}=\cj^W_{W_K}(6543322110+6534231201), K=D_7$;

$c=22$;  $2835_{22}=\cj^W_{W_K}(\e),    
K=\left(\sm\bul&.&\bul&\bul&\bul&\bul&\bul&\bul\\{}&{}&.&{}&{}&{}&{}&{}\esm\right)$;
                       
$c=22$;   $3200_{22}=\cj^W_{W_K}(\e), K=D_5\T(A_1\T A_1)$;

$c=22$;   $6075_{22}=\cj^W_{W_K}(5543322110\bxt\e),K=D_6\T A_1$;

$c=23$;    $4536_{23}=\cj^W_{W_K}(\e),
K=\left(\sm\bul&\bul&\bul&\bul&.&.&\bul&\bul\\{}&{}&\bul&{}&{}&{}&{}&{}\esm\right)$;

$c=24$; $4200_{24}+3360_{25}=\cj^W_{W_K}(6544322110+6544231201)$

$=\cj^W_{W_{K'}}(4200_{24}+3360_{25})$, $K=D_8,K'=E_8$;

$c=25$;  $2800_{25}+2100_{28}=\cj^W_{W_K}(120_{25}+105_{26}), K=E_7$;

$c=26$;  $840_{26}=\cj^W_{W_K}(\e),K=D_5A_3$;

$c= 26$;  $4096_{26}+4096_{27}=\cj^W_{W_K}(4096_{26}+4096_{27})
=\cj^W_{W_{K'}}(120_{25}\bxt\e+105_{26}\bxt\e)$,

$K=E_8,K'=E_7\T A_1$;

$c=28$;  $700_{28}=\cj^W_{W_K}(\e)=\cj^W_{W_{K'}}(\e)$,
$K=\left(\sm.&.&\bul&\bul&\bul&\bul&\bul&\bul\\{}&{}&\bul&{}&{}&{}&{}&{}\esm\right)$,
$K'=\left(\sm.&\bul&\bul&\bul&\bul&\bul&\bul&\bul\\{}&{}&.&{}&{}&{}&{}&{}\esm\right)$;

$c=28$;  $2240_{28}+840_{31}=\cj^W_{W_K}(2240_{28}+840_{31}),K=E_8$;

$c=29$;  $1400_{29}=\cj^W_{W_K}(\e)=\cj^W_{W_{K'}}(5544332110)$,
$K=\left(\sm\bul&.&\bul&\bul&\bul&\bul&\bul&\bul\\{}&{}&\bul&{}&{}&{}&{}&{}\esm\right)$, $K'=D_8$;

$c=30$;  $2268_{30}+1296_{33}=\cj^W_{W_K}(56_{30}+21_{33}),K=E_7$;

$c= 31$;  $3240_{31}+972_{32}=\cj^W_{W_K}(765443322110+756443322110)$

$=\cj^W_{W_{K'}}(56_{30}\bxt\e+35_{31}\bxt\e)$, $K=D_8,K'=E_7\T A_1$;

$c=32$;  $1400_{32}+2\T1575_{34}+350_{38}=\cj^W_{W_K}(1400_{32}+2\T1575_{34}+350_{38}),K=E_8$;

$c=34$;  $1050_{34}=\cj^W_{W_K}(665543322110), K=D_8$;

$c=36$;  $175_{36}=\cj^W_{W_K}(\e),K=A_8$;

$c=36$; $525_{36}=\cj^W_{W_K}(\e)$, $K=E_6$;

$c=37$;  $1400_{37}+2\T1008_{39}+56_{49}=\cj^W_{W_K}(1400_{37}+2\T1008_{39}+56_{49}),K=E_8$;

$c=38$;  $1344_{38}=\cj^W_{W_K}(27_{37}\bxt\e), K=E_7\T A_1$;

$c=39$;  $448_{39}=\cj^W_{W_K}(\e),K=E_6\T A_2$;

$c=42$;  $700_{42}+300_{44}=\cj^W_{W_K}(700_{42}+300_{44}),K=E_8$;

$c= 43$;   $400_{43}=\cj^W_{W_K}(77655443322110), K=D_8$;

$c=46$;   $567_{46}=\cj^W_{W_K}(7_{46}), K=E_7$;

$c=47$;  $560_{47}=\cj^W_{W_K}(7_{46}\bxt\e)=\cj^W_{W_{K'}}(560_{47}), K=E_7\T A_1,K'=E_8$;

$c=52$;  $210_{52}+160_{55}=\cj^W_{W_K}(210_{52}+160_{55}), K=E_8$;

$c=56$;   $50_{56}=\cj^W_{W_K}(8776655443322110), K=D_8$;

$c=63$;  $112_{63}+28_{68}=\cj^W_{W_K}(112_{63}+28_{68}), K=E_8$;

$c=64$;   $84_{64}=\cj^W_{W_K}(\e),K=E_7\T A_1$;

$c=74$;  $35_{74}=\cj^W_{W_K}(35_{74}),K=E_8$;

$c=91$;  $8_{91}=\cj^W_{W_K}(8_{91}), K=E_8$;

$c=120$;  $1_{120}=\cj^W_{W_K}(\e),K=E_8$.

\subhead 2.3. Table for type $E_7$\endsubhead

$c=0$;  $1_0=\cj^W_{W_K}(\e),K=\emp$;

$c=1$; $7_1=\cj^W_{W_K}(\e), K=\left(\sm\bul&.&.&.&.&.&.\\{}&{}&{}&.&{}&{}&{}\esm\right)$;

$c=2$;  $27_2=\cj^W_{W_K}(\e), K=\left(\sm\bul&.&\bul&.&.&.&.\\{}&{}&{}&.&{}&{}&{}\esm\right)$;

$c=3$;  $56_3+21_6=\cj^W_{W_K}(3210+3201),K=D_4$;

$c=3$;  $21_3=\cj^W_{W_K}(\e), K=\left(\sm.&.&.&.&\bul&.&\bul\\{}&{}&{}&\bul&{}&{}&{}\esm\right)$;

$c=4$;  $35_4=\cj^W_{W_K}(\e), K=\left(\sm\bul&.&\bul&.&\bul&.&\bul\\{}&{}&{}&.&{}&{}&{}\esm\right)$;
                 
$c=4$;  $120_4+105_5=\cj^W_{W_K}(3210\bxt\e+3201\bxt\e), K=D_4\T A_1$;

$c=5$;  $189_5+15_7=\cj^W_{W_K}(3210\bxt\e+3120\bxt\e), K=D_4\T(A_1\T A_1)$;

$c=6$;   $105_6=\cj^W_{W_K}(\e), K=\left(\sm\bul&\bul&.&.&\bul&.&\bul\\{}&{}&{}&\bul&{}&{}&{}\esm\right)$;

$c=6$;   $168_6=\cj^W_{W_K}(\e), K=\left(\sm\bul&\bul&.&.&.&\bul&\bul\\{}&{}&{}&.&{}&{}&{}\esm\right)$;     

$c=6$;   $210_6=\cj^W_{W_K}(\e), K=\left(\sm\bul&\bul&\bul&.&.&.&.\\{}&{}&{}&.&{}&{}&{}\esm\right)$;     
                 
$c=7$;   $315_7+2\T280_9+35_{13}=\cj^W_{W_K}(80_7+2\T90_9+20_{10}), K=E_6$;

$c=7$;   $189_7=\cj^W_{W_K}(\e), K=\left(\sm\bul&\bul&\bul&.&.&.&.\\
{}&{}&{}&\bul&{}&{}&{}\esm\right)$;     

$c=8$;   $280_8=\cj^W_{W_K}(\e), K=\left(\sm\bul&\bul&\bul&.&\bul&.&\bul\\
{}&{}&{}&.&{}&{}&{}\esm\right)$;     

$c=8$;   $405_8+189_{10}=\cj^W_{W_K}(432110\e+431201\e), K=D_5\T A_1$;

$c=9$;    $70_9=\cj^W_{W_K}(\e), K=A_2\T A_2\T A_2$;

$c=9$;   $216_9=\cj^W_{W_K}(\e), K=A_3\T A_1\T A_1$;

$c=9$;  $378_9=\cj^W_{W_K}(\e),K=\left(\sm\bul&\bul&\bul&.&.&\bul&\bul\\
{}&{}&{}&.&{}&{}&{}\esm\right)$;     
                     
$c=10$;  $420_{10}+336_{11}=\cj^W_{W_K}(432210+432201), K=D_6$;

$c= 10$; $210_{10}=\cj^W_{W_K}(\e)=\cj^W_{W_{K'}}(\e)$,
$K=\left(\sm\bul&\bul&\bul&.&.&\bul&\bul\\{}&{}&{}&\bul&{}&{}&{}\esm\right)$;

$c=11$;    $512_{11}+512_{12}=\cj^W_{W_K}(512_{11}+512_{12})$

$=\cj^W_{W_{K'}}(432210\bxt\e+432210\bxt\e),K=E_7, K'=D_6\T A_1$;

$c=12$;  $84_{12}=\cj^W_{W_K}(\e),K=A_3\T A_3$;

$c=12$;  $105_{12}=\cj^W_{W_K}(\e), K=D_4$;

$c=13$;  $210_{13}=\cj^W_{W_K}(\e)=\cj^W_{W_{K'}}(\e), K=A_3\T A_3\T A_1,K'=A_4\T A_2$;

$c=13$;  $420_{13}+336_{14}=\cj^W_{W_K}(54322110+54231201), K=D_6$;

$c=14$;  $378_{14}+84_{15}=\cj^W_{W_K}(54322110\bxt\e+53422110\bxt\e), K=D_6\T A_1$;

$c=15$;  $105_{15}=\cj^W_{W_K}(\e), K=\left(\sm\bul&\bul&\bul&\bul&.&.&.\\
{}&{}&{}&\bul&{}&{}&{}\esm\right)$;

$c=15$;  $405_{15}+189_{17}=\cj^W_{W_K}(30_{15}+15_{17}), K=E_6$;

$c=16$;  $315_{16}+2\T280_{18}+35_{22}=\cj^W_{W_K}(315_{16}+\T280_{18}+35_{22}),K=E_7$;

$c=16$;    $216_{16}=\cj^W_{W_K}(\e), K=\left(\sm\bul&\bul&\bul&\bul&\bul&.&\bul\\
{}&{}&{}&.&{}&{}&{}\esm\right)$;
                           
$c=17$;  $280_{17}=\cj^W_{W_K}(44332110\bxt\e), K=D_6\T A_1$;

$c=18$;  $70_{18}=\cj^W_{W_K}(\e),K=A_5\T A_2$;

$c=20$;  $189_{20}=\cj^W_{W_K}(\e), K=D_5$;

$c=21$;   $105_{21}=\cj^W_{W_K}(\e),K=A_6$;

$c=21$;   $168_{21}=\cj^W_{W_K}(\e), K=D_5\T A_1$;

$c=21$;  $210_{21}=\cj^W_{W_K}(5543322110), K=D_6$;

$c=22$;  $189_{22}=\cj^W_{W_K}(5543322110\bxt\e)=\cj^W_{W_{K'}}(189_{22}), K=D_6\T A_1,K'=E_7$;

$c=25$;  $120_{25}+105_{26}=\cj^W_{W_K}(120_{25}+105_{26}),K=E_7$;

$c=28$;  $15_{28}=\cj^W_{W_K}(\e), K=A_7$;

$c=30$;  $56_{30}+21_{33}=\cj^W_{W_K}(56_{30}+21_{33}), K=E_7$;

$c=31$;  $35_{31}=\cj^W_{W_K}(\e), K=D_6\T A_1$;

$c=36$;   $21_{36}=\cj^W_{W_K}(\e), K=E_6$;

$c=37$;   $27_{37}=\cj^W_{W_K}(27_{37}), K=E_7$;

$c=46$;   $7_{46}=\cj^W_{W_K}(7_{46}), K=E_7$;

$c=63$; $1_{63}=\cj^W_{W_K}(1_{63}), K=E_7$.

\subhead 2.4. Table for type $E_6$\endsubhead

$c=0$;  $1_0=\cj^W_{W_K}(\e), K=\emp$;

$c=1$;  $6_1=\cj^W_{W_K}(\e), K=\left(\sm \bul&.&.&.&.\\{}&{}&.&{}&{}\\{}&{}&.&{}&{}\esm\right)$;

$c=2$; $20_2=\cj^W_{W_K}(\e), K=\left(\sm \bul&.&\bul&.&.\\{}&{}&.&{}&{}\\{}&{}&.&{}&{}\esm\right)$;

$c=3$;  $30_3+15_5=\cj^W_{W_K}(3210+3201), K=D_4$;

$c=4$;   $15_4=\cj^W_{W_K}(\e), 
K=\left(\sm \bul&.&\bul&.&\bul\\{}&{}&.&{}&{}\\{}&{}&\bul&{}&{}\esm\right)$;

$c=4$;  $64_4=\cj^W_{W_K}(\e), K=\left(\sm \bul&\bul&.&\bul&.\\{}&{}&.&{}&{}\\{}&{}&.&{}&{}\esm\right)$;

$c=5$; $ 60_5=\cj^W_{W_K}(\e), K=\left(\sm \bul&\bul&.&\bul&.\\{}&{}&\bul&{}&{}\\{}&{}&.&{}&{}\esm\right)$;

$c=6$;  $24_6=\cj^W_{W_K}(\e), K=\left(\sm \bul&\bul&.&\bul&\bul\\{}&{}&.&{}&{}\\{}&{}&.&{}&{}\esm\right)$;

$c=6$;  $81_6=\cj^W_{W_K}(\e), K=\left(\sm \bul&\bul&\bul&.&.\\{}&{}&.&{}&{}\\{}&{}&.&{}&{}\esm\right)$;

$c=7$;  $80_7+2\T90_8+20_{10}=\cj^W_{W_K}(80_7+2\T90_8+20_{10}), K=E_6$;

$c=8$;  $60_8=\cj^W_{W_K}(\e), K=\left(\sm \bul&\bul&\bul&.&\bul\\{}&{}&.&{}&{}\\{}&{}&\bul&{}&{}\esm\right)$;

$c=9$;  $10_9=\cj^W_{W_K}(\e), K=A_2\T A_2\T A_2$;

$c=10$;   $81_{10}=\cj^W_{W_K}(\e), K=\left(\sm \bul&\bul&\bul&\bul&.\\{}&{}&.&{}&{}\\{}&{}&.&{}&{}\esm\right)$;

$c=11$; $60_{11}=\cj^W_{W_K}(\e), K=\left(\sm \bul&\bul&\bul&\bul&.\\{}&{}&.&{}&{}\\{}&{}&\bul&{}&{}\esm\right)$;

$c=12$;  $24_{12}=\cj^W_{W_K}(\e), K=D_4$;

$c=13$;   $64_{13}=\cj^W_{W_K}(44322110), K=D_5$;

$c=15$;  $30_{15}+15_{17}=\cj^W_{W_K}(30_{15}+15_{17}), K=E_6$;

$c= 16$;  $15_{16}=\cj^W_{W_K}(\e), K=A_5\T A_1$;

$c=20$;  $20_{20}=\cj^W_{W_K}(\e), K=D_5$;

$c=25$;  $6_{25}=\cj^W_{W_K}(6_{25}), K=E_6$;

$c=36$;   $1_{36}=\cj^W_{W_K}(1_{36}), K=E_6$;

\subhead 2.5. Table for type $F_4$\endsubhead

In this table we denote by $A_2$ (resp. $A'_2$) a subset $K\subsetneqq S_{af}$ of type $A_2$
which is contained (resp. not contained) in a subset of $S_{af}$ of type $A_3$.

$c=0$; $1_0=\cj^W_{W_K}(\e), K=\emp$;

$c=1$;  $4_1+2''_4=\cj^W_{W_K}(210+201), K=B_2$;           

$c=2$;  $9_2+2'_4=\cj^W_{W_K}(210\bxt\e+120\bxt\e), K=B_2\T A_1$;
 
$c=3$;   $8''_3=\cj^W_{W_K}(\e), K=A'_2$;

$c=3$;   $8'_3=\cj^W_{W_K}(\e)=\cj^W_{W_{K'}}(\e), K=A_2, K'=A_1\T A_1\T A_1$;

$c=4$;  $12_4+3\T9''_6+2\T6''_6+3\T1''_{12}=\cj^W_{W_K}(12_4+3\T9''_6+2\T6''_6+3\T1''_{12}), K=F_4$;

$c=5$;  $16_5+4''_7=\cj^W_{W_K}(32110\bxt\e+31201\bxt\e), K=C_3\T A_1$;

$c=6$;   $9'_6+4_8=\cj^W_{W_K}(32210+23120), K=B_4$;

$c=6$;  $6'_6=\cj^W_{W_K}(\e), K=A_2\T A'_2$;

$c=7$;  $4'_7=\cj^W_{W_K}(\e), K=A_3\T A_1$;

$c=9$;  $8'_9+1_{12}=\cj^W_{W_K}(4322110+4231201), K=B_4$;

$c=9$;  $8''_9=\cj^W_{W_K}(\e), K=C_3$;

$c=10$;  $9_{10}=\cj^W_{W_K}(\e)=\cj^W_{W_{K'}}(9_{10}), K=C_3\T A_1,K'=F_4$;

$c=13$;  $4_{13}+2''_{16}=\cj^W_{W_K}(4_{13}+2''_{16}), K=F_4$;

$c=16$;  $2_{16}=\cj^W_{W_K}(\e), K=B_4$;

$c=24$; $1_{24}=\cj^W_{W_K}(1_{24}), K=F_4$.

\subhead 2.6. Table for type $G_2$\endsubhead

$c=0$;  $1_0=\cj^W_{W_K}(\e), K=\emp$;

$c=1$;  $2_1+2\T1_3=\cj^W_{W_K}(2_1+2\T1_3), K=G_2$;

$c=2$;  $2_2=\cj^W_{W_K}(\e), K=A_1\T A_1$;

$c=3$;  $\ti 1_3=\cj^W_{W_K}(\e), K=A_2$;

$c=6$; $1_6=\cj^W_{W_K}(1_6), K=G_2$.

\head 3. Comments\endhead
\subhead 3.1\endsubhead
The tables in \S2 can be obtained by combining the induction/restriction tables in \cite{A05} with the explicit knowledge
of the $a$-function in \cite{L84}. This is how I first tried to obtain them. (While doing that I found that,
in \cite{L84}, the entry $[189_c]$ on p.364, line -5, should be replaced by $[189'_c]$.)
But eventually I used instead the induction/restriction tables in the CHEVIE package \cite{Ch}
combined with a package to compute the $a$-function
kindly supplied to me by Meinolf Geck. I am very grateful to Gongqin Li for doing the necessary programming.

\subhead 3.2\endsubhead
Now Propositions 1.4, 1.9 (and hence 1.6) follow (in the case of exceptional groups) from the tables in \S2
and the explicit determination of the Springer representations (see \cite{S85} and the references there).
The case of classical groups requires additional arguments; it will be considered elsewhere. One ingredient
in the proof is the interpretation \cite{L89} of the function $\g:\Irr W@>>>\NN$ in 1.1 in terms of the
$a$-function on $W_{af}$; this can be used to prove the inequality $a_{E'}\le \g(E)$ for any $E\in\Irr W$
and any $(K,E')\in\G(E)$. In addition to this, one has to use combinatorial arguments similar to those used in
\cite{L09}.

\subhead 3.3\endsubhead
Recall that in \cite{L86} an algorithm is given which determines the (ordinary) Green functions in terms
of the function $\g:\Irr W@>>>\NN$. Since by Proposition 1.4, this function is equal to the explicitly
computable function $c:\Irr W@>>>\NN$, we see that the Green functions can now be
calculated by a computer without the input of the Springer correspondence.

\subhead 3.4\endsubhead
We now assume that the characteristic $p$ of the ground field $\kk$
is a bad prime for $G$. Let $T$ be a maximal torus of $G$ and let $V=\QQ\ot\Hom(T,\kk^*)$,
$V^*=\QQ\ot\Hom(\kk^*,T)$. Let $R\sub V$ (resp. $\che R\sub V^*$) be  the set of roots (resp. coroots)
of $G$. For any $(e,e')\in V\T V^*$ let 
$W_{e,e'}$ be the reflection subgroup of the Weyl group $W$ (viewed as
a subgroup of $GL(V)$ or of $GL(V^*)$) defined in \cite{L15, 1.1} in terms of $R,\che R$.
Let $\cx$ be the set of reflection subgroups of $W$ of the form
$W_{e,e'}$ for some $(e,e')\in V\T V^*$ which is isolated in the sense of \cite{L15, 1.1}.
We define a subset $\cx^p$ of $\cx$ as follows. If $G$ is not of type $E_8$ we have $\cx^p=\cx$.
If $G$ is of type $E_8$ and $p=2$, $\cx^p$ consists of the subgroups in $\cx$ which are not
$W$-conjugate to one of type $A_2A_2A_2A_2$.
If $G$ is of type $E_8$ and $p=3$, $\cx^p$ consists of the subgroups in $\cx$ which are not
$W$-conjugate to one of type $D_6A_2,D_4D_4$ or $A_3A_3A_1A_1$.
It is known \cite{L15} that the irreducible representations of $W$ attached by the Springer correspondence to
unipotent classes in $G$ with the local system $\bbq$ are exactly those obtained by truncated
induction from special representations of the various subgroups in $\cx^p$.
For $E\in\Irr W$ let $\G^p(E)$ be the set of all pairs $(U,E')$ where $U\in\cx^p$ and $E'\in\Irr U$
is such that the multiplicity $(E:\Ind_U^W(E'))$ is nonzero. We have $\G^p(E)\ne\emp$. Let
$c^p_E=\max\{a_{E'};(U,E')\in\G^p(E)\}$; here $a_{E'}$ is computed in terms of the Weyl group $U$. Let
$\G^p_M(E)=\{(U,E')\in\G(E);a_{E'}=c^p_E\}$. We have $\G^p_M(E)\ne\emp$.
For $E\in \Irr W$ such that the unipotent class attached to $E$ under the Springer correspondence
is $\fc$, we define $\g^p(E)$ as the dimension of the Springer fibre at an element of $\fc$.

Experiments suggest that the following extension of 1.4(a) holds.

(a) For any $E\in\Irr W$ we have $c^p_E=\g^p(E)$.

\enddocument